\documentclass[14pt]{article}
\usepackage{amsfonts}
\usepackage{amsmath}
\usepackage{graphicx}

\newtheorem{theorem}{Theorem}

\newtheorem{proposition}[theorem]{Proposition}

\begin{document}

\title{ \textbf{Rent seeking games with tax evasion }}
\author{Olivia Bund\u{a}u$^a$, Mihaela Neam\c{t}u$^b$,
Dumitru Opri\c{s}$^c$}

\date { }
\maketitle

\begin{tabular}{cccccccc}
\scriptsize{$^{a}$ Department of Mathematics, "Politehnica"
University of Timi\c soara, E-mail: obundau@yahoo.es,}\\
%\scriptsize{West University of Timi\c soara, Bd. V. Parvan, nr. 4, 300223, Timi\c soara, Romania,}\\
%\scriptsize{E-mail: obundau@yahoo.es,}\\
\scriptsize{$^{b}$Department of Economic Informatics and Statistics, Faculty of Economics,}\\
\scriptsize{West University of Timi\c soara, e-mail:mihaela.neamtu@fse.uvt.ro,}\\
%\scriptsize{E-mail:mihaela.neamtu@fse.uvt.ro,}\\
\scriptsize{$^{c}$ Department of Applied
Mathematics, Faculty of Mathematics,}\\
\scriptsize{West University of Timi\c soara, e-mail: opris@math.uvt.ro.}\\
%\scriptsize{E-mail: opris@math.uvt.ro}\\

\end{tabular}

\begin{abstract}
We consider the static and dynamic models of Cournot duopoly with tax
evasion. In the dynamic model we introduce the time delay and we analyze the
local stability of the stationary state. There is a critical value of the
delay when the Hopf bifurcation occurs.
\end{abstract}

\textit{{\small Mathematics Subject Classification:\, 34K18,
47N10; Jel Classification: C61, C62, H26 }
}

\section{\protect\bigskip Introduction}

\hspace{0.4cm}During the last decades revenues from indirect tax
have become increasingly important in many economies. Substantial
attention has been devoted to evasion of indirect taxes. It is
well known that indirect tax evasion, especially evasion of VAT,
may erode a substantial part of tax revenues [2], [4], [5].

In [3] a model with tax evasion is presented. The authors consider
n firms which enter the market with a homogenous good. These firms
have to pay an ad valorem sales tax, but may evade a certain
amount of their tax duty. The aims of the firms are to maximize
their profits. The equilibrium point is determined and an economic
analysis is made.

Based on [1], [3], [7], [8], [10], in our paper we present three
economic models with tax evasion: the static model of Cournot
duopoly with tax evasion in Section 2, the dynamic model of
Cournot duopoly with tax evasion in Section 3 and the dynamic
model with tax evasion and time delay in Section 4.

In Section 2, in the static model the purpose of the firms is to
maximize their profits. We determine the firms' outputs and the
declared revenues which maximize the profits, as well as the
conditions for the model's parameters in which the maxim profits
are obtained. Using Maple 11, the variables orbits are displayed.

In Section 3, the dynamic model describes the variation in time of
the firms' outputs and the declared revenues. We study the local
stability for the stationary state and the conditions under which
it is asymptotically stable.

In Section 4, we formulate a new dynamic model, based on the model
from Section 3, in which the time delay is introduced. That means,
the two firms do not enter the market at the same time. One of
them is the leader firm and the other is the follower firm. The
second one knows the leader's output in the previous moment
$t-\tau ,$ $\tau \geq 0.$

Using classical methods [6], [9] we investigate the local
stability of the stationary state by analyzing the corresponding
transcendental characteristic equation of the linearized system.
By choosing the delay as a bifurcation parameter we show that this
model exhibits a limit cycle.

Finally numerical simulations, some conclusions and future
research possibilities are offered.

\section{The static model of Cournot duopoly with tax evasion}

\qquad The static model of Cournot duopoly is described by an economic game,
where two firms enter the market with a homogenous consumption product. The
elements which describe the model are: the quantities which enter the market
from the two firms $x_{i}\geq 0,$ $i=\overline{1,2};$ the declared revenues $%
z_{i},$ $i=\overline{1,2}$; the inverse demand function $p:\mathbb{R}%
_{+}\rightarrow \mathbb{R}_{+}$ ( $p$ is a derivable function with $%
p^{\prime }\left( x\right) <0,\underset{x\rightarrow a_{1}}{lim}p\left(
x\right) =0,$ $\underset{x\rightarrow 0}{lim}p\left( x\right) =b_{1},\left(
a_{1}\in \overline{\mathbb{R}},b_{1}\in \overline{\mathbb{R}}\right) $; the
penalty function $F:\mathbb{R}_{+}\rightarrow \mathbb{R}_{+}$ ($F$ is a
derivable function with $F^{\prime }\left( x\right) >0,$ $F^{\prime \prime
}\left( x\right) >0,$ $F\left( 0\right) =0);$ the cost functions $C_{i}:%
\mathbb{R}_{+}\rightarrow \mathbb{R}_{+}$ ( $C_{i}$ are derivable functions
with $C_{i}^{\prime }\left( x_{i}\right) >0,$ $C_{i}^{\prime \prime }\geq 0,$
$i=\overline{1,2}$ ). The government levies an ad valorem tax on each firm's
sales at the rate $t_{1}\in \left( 0,1\right) $ and $q\in \left[ 0,1\right] $
is the probability with which the tax evasion is detected.

The true tax base of firm $i$ is $x_{i}p\left( x_{1}+x_{2}\right) .$ Firm $i$
declares $z_{i}\leq x_{i}p\left( x_{1}+x_{2}\right) $ as tax base to the tax
authority. Accordingly, evaded revenues of firm $i$ are given by $%
x_{i}p\left( x_{1}+x_{2}\right) -z_{i}.$ With probability $1-q$ tax evasion
remains undetected and the tax bill of firm $i$ amounts to $t_{1}z_{i}.$ The
tax authority detects tax evasion of firm $i$ with probability $q$. In case
of detection, firm $i$ has to pay taxes on the full amount of revenues, $%
x_{i}p\left( x_{1}+x_{2}\right) ,$ and, in addition, a penalty
$F\left( x_{i}p\left( x_{1}+x_{2}\right) -z_{i}\right) $. The
penalty is increasing and convex in evaded revenues $x_{i}p\left(
x_{1}+x_{2}\right) -z_{i}$. Moreover, it is assumed that $F\left(
0\right) =0$, namely law-abiding firms go unpunished.

The profit functions of the two firms are: $P_{i}:\mathbb{R}%
_{+}^{2}\rightarrow \mathbb{R}_{+},$ $i=\overline{1,2}$, given by:%
\begin{eqnarray}
P_{i} &=&P_{i}\left( x_{1},x_{2},z_{1},z_{2}\right) =\left( 1-q\right) \left[
x_{i}p\left( x_{1}+x_{2}\right) -C_{i}\left( x_{i}\right) -t_{1}z_{i}\right]
+  \notag \\
&&q\left[ \left( 1-t_{1}\right) x_{i}p\left( x_{1}+x_{2}\right)
-C_{i}\left( x_{i}\right) -F\left( x_{i}p\left( x_{1}+x_{2}\right)
-z_{i}\right) \right]. \quad
\end{eqnarray}%
The first bracketed term in $\left( 1\right) $ equals the profit of firm $i$
if evasion activities remain undetected. The second term in $\left( 1\right)
$ represents the profit of firm $i$ in case tax evasion is detected.

The firm's aim is to maximize $\left( 1\right) $ with respect to output $%
x_{i}$ and declared revenues $z_{i}.$ This aim represents a mathematical
optimization problem which is given by:%
\begin{equation}
\underset{\{x_{i},z_{i}\}}{max}P_{i},\quad i=\overline{1,2}.
\end{equation}%
From the hypothesis about the functions $p,F,C_{i},i=\overline{1,2}$ , we
have:

\begin{proposition}
The solution of problem $\left( 2\right) $ is given by the solution of the
following system:%
\begin{eqnarray}
\dfrac{\partial P_{i}}{\partial x_{i}} &=&\left[ 1-qt_{1}-qF^{\prime }\left(
x_{i}p\left( x_{1}+x_{2}\right) -z_{i}\right) \right] \left[ p\left(
x_{1}+x_{2}\right) +x_{i}p^{\prime }\left( x_{1}+x_{2}\right) \right]
-C_{i}\left( x_{i}\right) =0  \notag \\
\dfrac{\partial P_{i}}{\partial z_{i}} &=&-\left( 1-q\right)
t_{1}+qF^{\prime }\left( x_{i}p\left( x_{1}+x_{2}\right)
-z_{i}\right) =0,\quad i=\overline{1,2}.
\end{eqnarray}
\end{proposition}

In what follows, we will consider the penalty function $F\left( x\right) =%
\dfrac{1}{2}st_{1}x^{2},$ $s\geq 1,$ the cost functions $C_{i}$
$\left(
x_{i}\right) =c_{i}x_{i},$ $c_{i}>0,$ $i=1,2$ and the price function $%
p\left( x\right) =\dfrac{1}{x}$ .

From $\left( 3\right) $ we can deduce:

\begin{proposition}
If $\dfrac{1-q}{qs+q-1}c_{1}\leq c_{2}\leq \dfrac{qs+q-1}{1-q}c_{1},$ then
the solution of system $\left( 3\right) $ is given by :%
\begin{eqnarray}
x_{1}^{\ast } &=&\dfrac{c_{2}\left( 1-t_{1}\right) }{\left(
c_{1}+c_{2}\right) ^{2}},\qquad \qquad \quad x_{2}^{\ast }=\dfrac{%
c_{1}\left( 1-t_{1}\right) }{\left( c_{1}+c_{2}\right) ^{2}},  \notag \\
z_{1}^{\ast } &=&\dfrac{c_{2}}{c_{1}+c_{2}}-\dfrac{1-q}{qs},\qquad
z_{2}^{\ast }=\dfrac{c_{1}}{c_{1}+c_{2}}-\dfrac{1-q}{qs}.
\end{eqnarray}
\end{proposition}

For the parameters $c_1=0.3$, $c_2=0.6$, $q=0.12$, $t_1=0.16$, the
variations of the variables $z_1$, $z_2$, and the profits $P_1$,
$P_2$ are given in the following figures:

\begin{center}\begin{tabular}{ccc}
%\hline
\includegraphics[width=4cm]{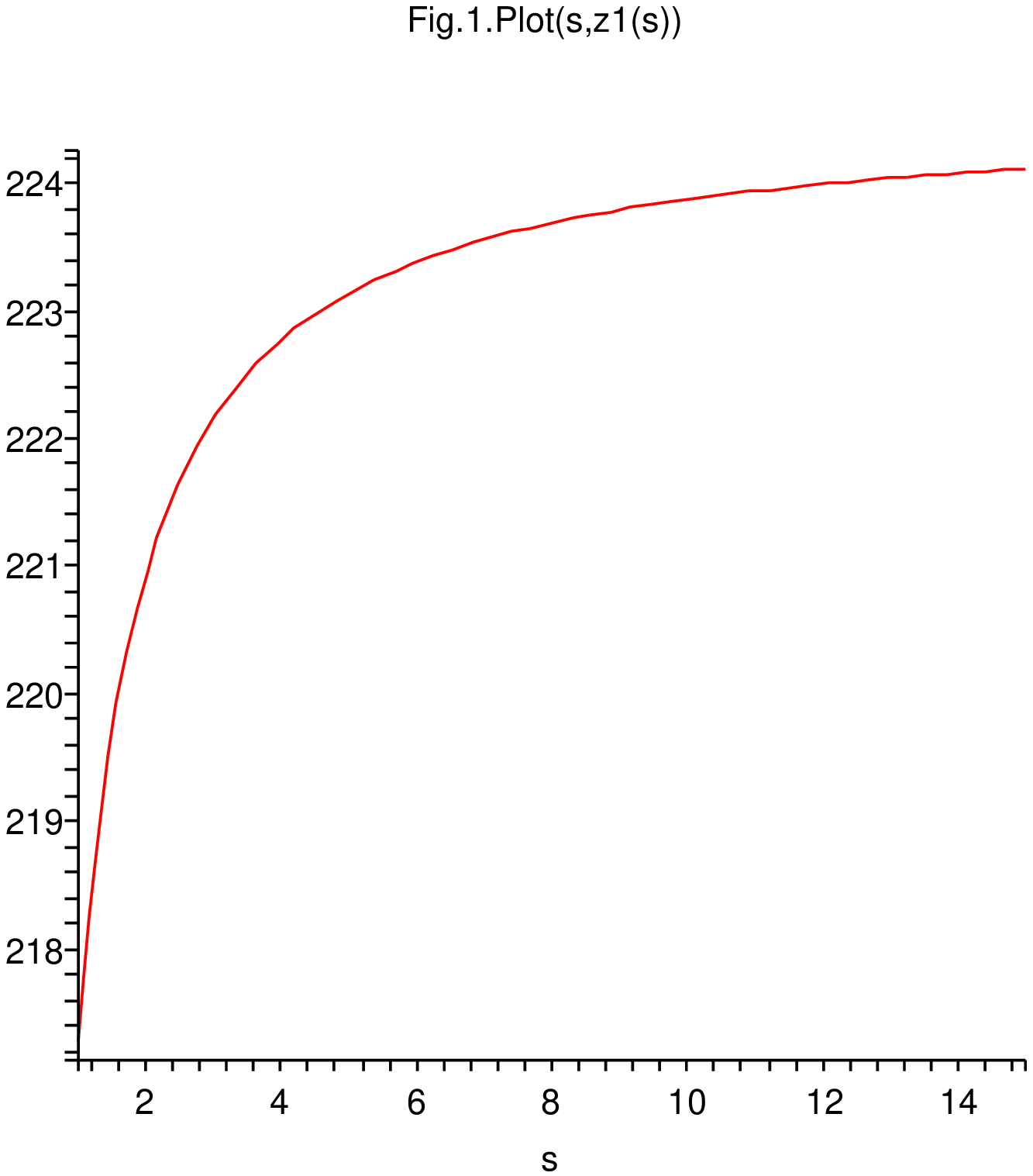}
\includegraphics[width=4cm]{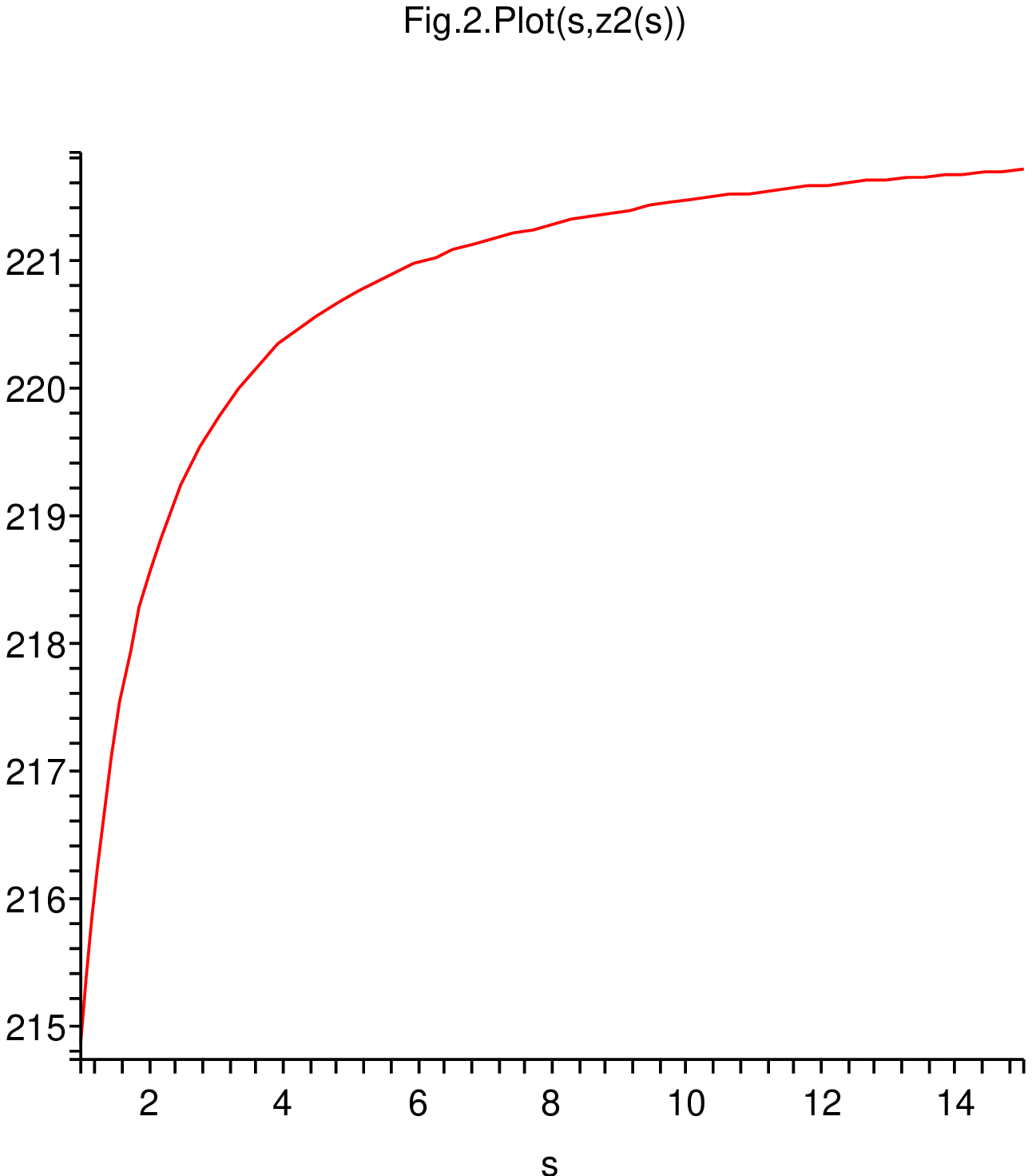}

\end{tabular}
\end{center}

\begin{center}\begin{tabular}{ccc}
%\hline
\includegraphics[width=5cm]{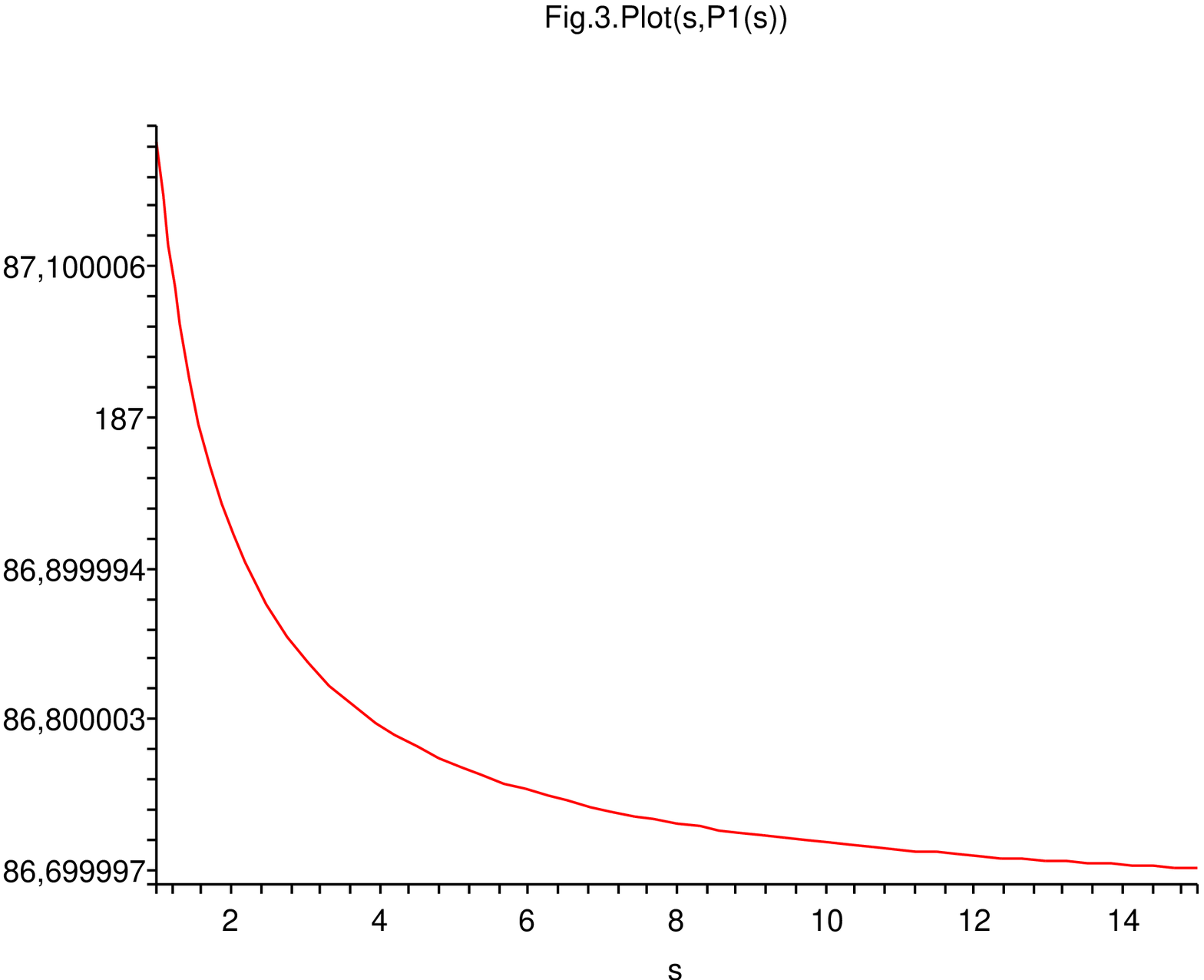}
\includegraphics[width=5cm]{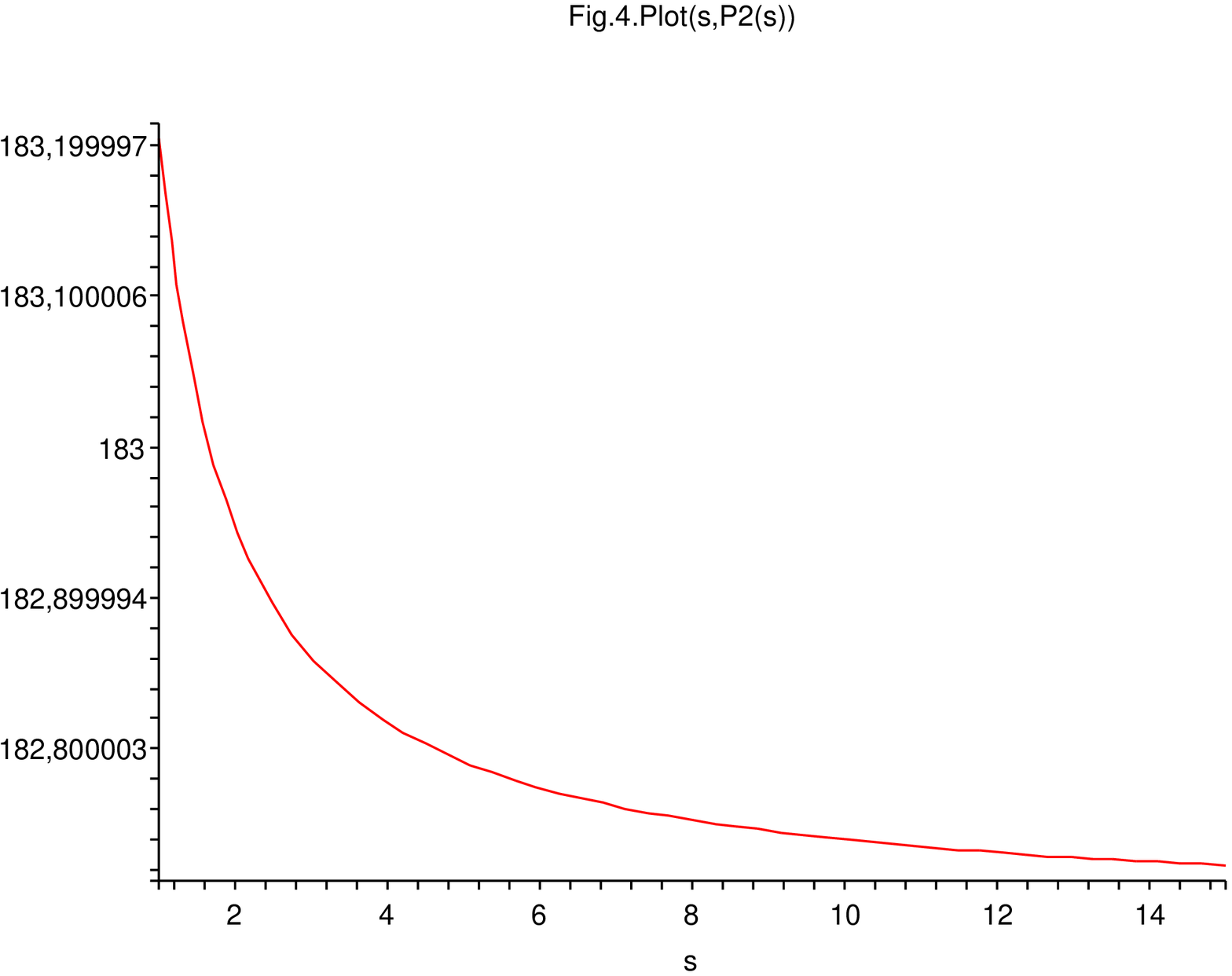}

\end{tabular}
\end{center}

\section{The dynamic model of Cournot duopoly with tax evasion}

The dynamic model describes the variation in time of output $x_{i}\left(
t\right) ,$ $i=1,2$ taking into account the marginal profits $\dfrac{%
\partial P_{i}}{\partial x_{i}},$ $i=1,2.$ Assume that each agent adjusts
its declared revenue $z_{i}\left( t\right) $, $i=1,2$
proportionally to the marginal profits $\dfrac{\partial
P_{i}}{\partial z_{i}},$ $i=1,2.$ Then,
the dynamic model is given by the following differential system of equations:%
\begin{eqnarray}
\overset{\cdot }{x}_{i}\left( t\right) &=&k_{i}\dfrac{\partial P_{i}}{%
\partial x_{i}}=k_{i}\{\left[ 1-qt_{1}-qF^{\prime }\left( x_{i}p\left(
x_{1}+x_{2}\right) -z_{i}\right) \right] \cdot  \notag \\
&&\left[ p\left( x_{1}+x_{2}\right) +x_{i}p^{\prime }\left(
x_{1}+x_{2}\right) \right] -C_{i}\left( x_{i}\right)\} ,\quad \\
\overset{\cdot }{z}_{i}\left( t\right) &=&h_{i}\dfrac{\partial P_{i}}{%
\partial z_{i}}=h_{i}\left[ -\left( 1-q\right) t_{1}+qF^{\prime }\left(
x_{i}p\left( x_{1}+x_{2}\right) -z_{i}\right) \right] ,\quad i=\overline{1,2}
\notag
\end{eqnarray}%
with the initial conditions $x_{i}\left( 0\right) =x_{i0},$ $z_{i}\left(
0\right) =z_{i0},$ $i=\overline{1,2}$ and $h_{i}>0,$ $k_{i}>0,$ $i=\overline{%
1,2}.$

For $F\left( x\right) =\dfrac{1}{2}st_{1}x^{2},$ $s\geq 1$ and $C_{i}\left(
x_{i}\right) =c_{i}x_{i},$ $c_{i}>0,$ $i=\overline{1,2}$ system $\left(
5\right) $ becomes :

\begin{eqnarray}
\overset{\cdot }{x}_{1}\left( t\right) &=&k_{1}\{\left[ 1-qt_{1}-qst_{1}%
\left( x_{1}\left( t\right) p\left( x_{1}\left( t\right) +x_{2}\left(
t\right) \right) -z_{1}\left( t\right) \right) \right] \cdot  \notag \\
&&\cdot \left[ p\left( x_{1}\left( t\right) +x_{2}\left( t\right) \right)
+x_{1}\left( t\right) p^{\prime }\left( x_{1}\left( t\right) +x_{2}\left(
t\right) \right) \right] -c_{1}\}  \notag \\
\overset{\cdot }{x}_{2}\left( t\right) &=&k_{2}\{\left[ 1-qt_{1}-qst_{1}%
\left( x_{2}\left( t\right) p\left( x_{1}\left( t\right) +x_{2}\left(
t\right) \right) -z_{2}\left( t\right) \right) \right] \cdot  \notag \\
&&\cdot \left[ p\left( x_{1}\left( t\right) +x_{2}\left( t\right) \right)
+x_{2}\left( t\right) p^{\prime }\left( x_{1}\left( t\right) +x_{2}\left(
t\right) \right) \right] -c_{2}\} \\
\overset{\cdot }{z}_{1}\left( t\right) &=&h_{1}\left[ -\left( 1-q\right)
t_{1}+qst_{1}\left( x_{1}\left( t\right) p\left( x_{1}\left( t\right)
+x_{2}\left( t\right) \right) -z_{1}\left( t\right) \right) \right]  \notag
\\
\overset{\cdot }{z}_{2}\left( t\right) &=&h_{2}\left[ -\left( 1-q\right)
t_{1}+qst_{1}\left( x_{2}\left( t\right) p\left( x_{1}\left( t\right)
+x_{2}\left( t\right) \right) -z_{2}\left( t\right) \right) \right]  \notag
\\
x_{i}\left( 0\right) &=&x_{i0},\text{ }z_{i}\left( 0\right) =z_{i0},\text{ }%
i=1,2.  \notag
\end{eqnarray}

System $\left( 6\right) $ has the stationary state $\left(
x_{1}^{\ast},x_{2}^{\ast },z_{1}^{\ast },z_{2}^{\ast }\right) $
given by $\left( 4\right) .$

Let $u_{1}\left( t\right) =x_{1}\left( t\right) -x_{1}^{\ast }$,
$u_{2}\left( t\right) =x_{2}\left( t\right) -x_{2}^{\ast }$,
$u_{3}\left( t\right) =z_{1}\left( t\right) -z_{1}^{\ast }$,
$u_{4}\left( t\right) =z_{2}\left( t\right) -z_{2}^{\ast }.$

By expanding $\left( 6\right) $ in a Taylor series around the stationary
state $\left( x_{1}^{\ast },x_{2}^{\ast },z_{1}^{\ast },z_{2}^{\ast }\right)
$ and neglecting the terms of higher order than the first order, we have the
following linear appro\-xi\-ma\-tion of system $\left( 6\right) $ :%
\begin{eqnarray}
\overset{\cdot }{u}_{1}\left( t\right) &=&k_{1}\{a_{10}u_{1}\left( t\right)
+a_{01}u_{2}\left( t\right) +a_{001}u_{3}\left( t\right) \}  \notag \\
\overset{\cdot }{u}_{2}\left( t\right) &=&k_{2}\{b_{10}u_{1}\left( t\right)
+b_{01}u_{2}\left( t\right) +b_{001}u_{4}\left( t\right) \}  \notag \\
\overset{\cdot }{u}_{3}\left( t\right) &=&h_{1}\{c_{10}u_{1}\left( t\right)
+c_{01}u_{2}\left( t\right) +c_{001}u_{3}\left( t\right) \} \\
\overset{\cdot }{u}_{4}\left( t\right) &=&h_{2}\{d_{10}u_{1}\left( t\right)
+d_{01}u_{2}\left( t\right) +d_{001}u_{4}\left( t\right) \}  \notag
\end{eqnarray}%
where:

\begin{eqnarray}
a_{10} &=&-qst_{1}\dfrac{c_{1}^{2}}{\left( 1-t_{1}\right)
^{2}}+(1-t_1)\left( 2p^{\prime }\left( x_{1}^{\ast }+x_{2}^{\ast
}\right) +x_{1}^{\ast }p^{\prime \prime }\left( x_{1}^{\ast
}+x_{2}^{\ast }\right) \right),  \notag
\\
a_{01} &=&-qst_{1}x_{1}^{\ast }p^{\prime }\left( x_{1}^{\ast
}+x_{2}^{\ast }\right) \dfrac{c_{1}}{1-t_{1}}+(1-t_1)\left(
p^{\prime }\left( x_{1}^{\ast }+x_{2}^{\ast }\right) +x_{1}^{\ast
}p^{\prime \prime }\left(
x_{1}^{\ast }+x_{2}^{\ast }\right) \right),  \notag \\
a_{001} &=&\dfrac{qst_{1}c_{1}}{1-t_{1}},\quad b_{001}=\dfrac{qst_{1}c_{2}}{%
1-t_{1}}, \quad c_{001}=d_{001}=-qst_{1}, \quad c_{10}=\dfrac{qst_{1}c_{1}}{1-t_{1}},\\
b_{10} &=&-qst_{1}x_{2}^{\ast }p^{\prime }\left( x_{1}^{\ast
}+x_{2}^{\ast }\right) \dfrac{c_{2}}{1-t_{1}}+(1-t_1)\left(
p^{\prime }\left( x_{1}^{\ast }+x_{2}^{\ast }\right) +x_{2}^{\ast
}p^{\prime \prime }\left(
x_{1}^{\ast }+x_{2}^{\ast }\right) \right),  \notag \\
b_{01} &=&-qst_{1}\dfrac{c_{2}^{2}}{\left( 1-t_{1}\right)
^{2}}+(1-t_1)\left( 2p^{\prime }\left( x_{1}^{\ast }+x_{2}^{\ast
}\right) +x_{2}^{\ast }p^{\prime \prime }\left( x_{1}^{\ast
}+x_{2}^{\ast }\right) \right),  \notag
\\
c_{01} &=&qst_{1}x_{1}^{\ast }p^{\prime }\left( x_{1}^{\ast
}+x_{2}^{\ast }\right), \quad d_{10}=qst_1x_{2}^{\ast }p^{\prime
}\left( x_{1}^{\ast }+x_{2}^{\ast }\right), \quad
d_{01}=qst_1\dfrac{c_{2}}{1-t_{1}}. \notag
\end{eqnarray}

The characteristic equation associated to $\left( 7\right) $ is given by:%
\begin{equation}
\lambda ^{4}+m_{43}\lambda ^{3}+m_{42}\lambda ^{2}+m_{41}\lambda +m_{40}=0
\end{equation}%
where:%
\begin{eqnarray}
m_{43} &=&-k_{1}a_{10}-k_{2}b_{01}-h_{1}c_{001}-h_{2}d_{001}  \notag \\
m_{42} &=&k_{1}k_{2}a_{10}b_{01}+\left( k_{1}a_{10}+k_{2}b_{01}\right)
\left( h_{1}c_{001}+h_{2}d_{001}\right) -h_{1}k_{1}a_{001}c_{10}-  \notag \\
&&-k_{2}h_{2}b_{001}d_{01}+h_{1}h_{2}c_{001}d_{001}-k_{1}k_{2}a_{01}b_{10}
\notag \\
m_{41} &=&k_{1}k_{2}a_{10}b_{01}\left(
h_{1}c_{001}-h_{2}d_{001}\right)
-k_{2}h_{1}h_{2}c_{001}d_{001}b_{01}+ \notag\\
&&+h_{1}k_{1}c_{10}a_{001}\left( k_{2}b_{01}+h_{2}d_{001}\right)
-k_{1}k_{2}h_{2}b_{001}a_{01}d_{10}+ \notag\\
&&+k_{1}k_{2}h_{2}a_{10}b_{001}d_{01}+k_{2}h_{1}h_{2}b_{001}c_{001}d_{01}+ \notag\\
&&+k_{1}k_{2}a_{01}b_{10}\left( h_{2}d_{001}+h_{1}c_{001}\right)-k_{1}k_{2}h_{1}a_{001}b_{10}c_{01} \\
m_{40}
&=&k_{1}k_{2}h_{1}h_{2}(a_{10}b_{01}c_{001}d_{001}-a_{001}b_{01}c_{10}d_{001}+a_{001}b_{001}c_{10}d_{01}+
\notag \\
&&+a_{01}b_{001}c_{001}d_{10}-a_{001}b_{001}c_{01}d_{10}-a_{10}b_{001}c_{001}d_{01}+
\notag \\
&&+a_{001}b_{10}c_{01}d_{001}-a_{01}b_{10}c_{001}d_{001}).  \notag
\end{eqnarray}

A necessary and sufficient condition as equation $\left( 9\right) $ has all
roots with negative real part is given by Routh-Hurwitz criterion:%
%\begin{eqnarray}
%D_{1} &=&m_{43}>0,\text{ \qquad \qquad\ \qquad }  \notag \\
%D_{2} &=&m_{43}m_{42}-m_{41}>0,\text{ } \\
%D_{3} &=&m_{41}D_{2}-m_{43}^{2}m_{40}>0,\text{ }  \notag \\
%D_{4} &=&m_{41}D_{3}>0.  \notag
%\end{eqnarray}

\begin{eqnarray}
D_{1}=m_{43}>0,D_{2}=m_{43}m_{42}-m_{41}>0,D_{3}=m_{41}D_{2}-m_{43}^{2}m_{40}>0,D_{4}
=m_{41}D_{3}>0.
\end{eqnarray}

From the Routh-Hurwitz criterion we have:

\begin{proposition}
The stationary state of system $\left( 6\right) $ is
asymptotically stable if and only if conditions $\left( 11\right)
$ hold.
\end{proposition}

\section{The dynamic model with tax evasion and time delay}

\qquad In [7] and [8] we have studied the rent seeking games with
time delay and distributed delay. In the present section we
analyze the rent seeking games with tax evasion and delay. For
$\tau=0$ we obtain the model from [3]. For $\tau=0$ and $t_1=0$ we
obtain the model from [1]. We consider the model from section 3
where we introduce the time delay $\tau$.  We suppose the first
firm is the leader and the second firm is the follower. The
follower knows the quantity of the leader firm, $x_{1}\left(
t-\tau \right) ,$ which entered the market at the moment $t-\tau
,$ $\tau >0.$

The differential system which describes this model is given by:%
\begin{eqnarray}
\overset{\cdot }{x}_{1}\left( t\right) &=&k_{1}\{\left[ 1-qt_{1}-qst_{1}%
\left( x_{1}\left( t\right) p\left( x_{1}\left( t\right) +x_{2}\left(
t\right) \right) -z_{1}\left( t\right) \right) \right] \cdot  \notag \\
&&\left[ p\left( x_{1}\left( t\right) +x_{2}\left( t\right) \right)
+x_{1}\left( t\right) p^{\prime }\left( x_{1}\left( t\right) +x_{2}\left(
t\right) \right) \right] -c_{1}\}  \notag \\
\overset{\cdot }{x}_{2}\left( t\right) &=&k_{2}\{\left[ 1-qt_{1}-qst_{1}%
\left( x_{2}\left( t\right) p\left( x_{1}\left( t-\tau \right) +x_{2}\left(
t\right) \right) -z_{2}\left( t\right) \right) \right] \cdot  \notag \\
&&\left[ p\left( x_{1}\left( t-\tau \right) +x_{2}\left( t\right) \right)
+x_{2}\left( t\right) p^{\prime }\left( x_{1}\left( t-\tau \right)
+x_{2}\left( t\right) \right) \right] -c_{2}\} \\
\overset{\cdot }{z}_{1}\left( t\right) &=&h_{1}\left[ -\left( 1-q\right)
t_{1}+qst_{1}\left( x_{1}\left( t\right) p\left( x_{1}\left( t\right)
+x_{2}\left( t\right) \right) -z_{1}\left( t\right) \right) \right]  \notag
\\
\overset{\cdot }{z}_{2}\left( t\right) &=&h_{2}\left[ -\left( 1-q\right)
t_{1}+qst_{1}\left( x_{2}\left( t\right) p\left( x_{1}\left( t\right)
+x_{2}\left( t\right) \right) -z_{2}\left( t\right) \right) \right]  \notag
\\
x_{1}\left( \theta \right) &=&\varphi \left( \theta \right) ,\text{ }\theta
\in \left[ -\tau ,0\right] ,\text{ }x_{2}\left( 0\right) =x_{20},\text{ }%
z_{i}\left( 0\right) =z_{i0},\text{ }k_{i}>0,\text{
}h_{i}>0,\text{ }i=1,2. \notag
\end{eqnarray}%
For $p\left( x\right) =\dfrac{1}{x}$ the stationary state of system $\left(
12\right) $ is given by $\left( 4\right) .$

With respect to the transformation $u_{1}\left( t\right)
=x_{1}\left( t\right) -x_{1}^{\ast }$, $u_{2}\left( t\right)
=x_{2}\left( t\right) -x_{2}^{\ast }$, $u_{3}\left( t\right)
=z_{1}\left( t\right) -z_{1}^{\ast }$, $u_{4}\left( t\right)
=z_{2}\left( t\right) -z_{2}^{\ast },$ and by expanding $\left(
12\right) $ in a Taylor series around the stationary state $\left(
x_{1}^{\ast },x_{2}^{\ast },z_{1}^{\ast },z_{2}^{\ast }\right) $
and neglecting the terms of higher order than the first order, we
obtain the following linear approximation of system $\left(
12\right) $ :

\begin{eqnarray}
\overset{\cdot }{u}_{1}\left( t\right) &=&k_{1}\{a_{10}u_{1}\left( t\right)
+a_{01}u_{2}\left( t\right) +a_{001}u_{3}\left( t\right) \}  \notag \\
\overset{\cdot }{u}_{2}\left( t\right) &=&k_{2}\{b_{10}u_{1}\left( t-\tau
\right) +b_{01}u_{2}\left( t\right) +b_{001}u_{4}\left( t\right) \}  \notag
\\
\overset{\cdot }{u}_{3}\left( t\right) &=&h_{1}\{c_{10}u_{1}\left( t\right)
+c_{01}u_{2}\left( t\right) +c_{001}u_{3}\left( t\right) \} \\
\overset{\cdot }{u}_{4}\left( t\right) &=&h_{2}\{d_{10}u_{1}\left( t\right)
+d_{01}u_{2}\left( t\right) +d_{001}u_{4}\left( t\right) \}  \notag
\end{eqnarray}%
where $%
a_{10},a_{01},a_{001},b_{10},b_{01},b_{001},c_{10},c_{01},c_{001},d_{10},d_{01},d_{001},
$ are given by $\left( 8\right) .$

The corresponding characteristic equation of $\left( 13\right) $ is :%
\begin{equation}
\lambda ^{4}+n_{43}\lambda ^{3}+n_{42}\lambda ^{2}+n_{41}\lambda
+n_{40}+e^{-\lambda \tau }\left( n_{22}\lambda ^{2}+n_{21}\lambda
+n_{20}\right) =0,
\end{equation}%
where
\begin{eqnarray*}
n_{43} &=&m_{43},\text{ }n_{22}=-k_{1}k_{2}a_{01}b_{10} \\
n_{42} &=&k_{1}k_{2}a_{10}b_{01}+\left( k_{1}a_{10}+k_{2}b_{01}\right)
\left( h_{1}c_{001}+h_{2}d_{001}\right) -h_{1}k_{1}a_{001}c_{10}- \\
&&-k_{2}h_{2}b_{001}d_{01}+h_{1}h_{2}c_{001}d_{001} \\
n_{41} &=&k_{1}k_{2}a_{10}b_{01}\left(
h_{1}c_{001}-h_{2}d_{001}\right)
-k_{2}h_{1}h_{2}c_{001}d_{001}b_{01}+ \\
&&+h_{1}k_{1}c_{10}a_{001}\left( k_{2}b_{01}+h_{2}d_{001}\right)
-k_{1}k_{2}h_{2}b_{001}a_{01}d_{10}+ \\
&&+k_{1}k_{2}h_{2}a_{10}b_{001}d_{01}+k_{2}h_{1}h_{2}b_{001}c_{001}d_{01}\\
n_{40}
&=&k_{1}k_{2}h_{1}h_{2}(a_{10}b_{01}c_{001}d_{001}-a_{001}b_{01}c_{10}d_{001}+a_{001}b_{001}c_{10}d_{01}+
\\
&&+a_{01}b_{001}c_{001}d_{10}-a_{001}b_{001}c_{01}d_{10}-a_{10}b_{001}c_{001}d_{01})
\\
n_{21} &=&k_{1}k_{2}a_{01}b_{10}\left( h_{2}d_{001}+h_{1}c_{001}\right)-k_{1}k_{2}h_{1}a_{001}b_{10}c_{01} \\
n_{20}
&=&k_{1}k_{2}h_{1}h_{2}(a_{001}b_{10}c_{01}d_{001}-a_{01}b_{10}c_{001}d_{001}).
\end{eqnarray*}%

The roots of $\left( 14\right) $ depend on $\tau .$ Considering
$\tau $ as parameter, we determine $\tau _{0}$ so that $\lambda
=i\omega $ is a root of $\left( 14\right) .$ Substituting $\lambda
=i\omega $ into equation $\left( 14\right) $ we obtain:\newline
$\omega ^{4}-in_{43}\omega ^{3}-n_{42}\omega ^{2}+in_{41}\omega
+n_{40}+\left( -n_{22}\omega ^{2}+in_{21}\omega +n_{20}\right)
\left( cos\omega \tau -isin\omega \tau \right) =0.$

From the above equation we have:%
\begin{equation}
\omega ^{8}+r_{6}\omega ^{6}+r_{4}\omega ^{4}+r_{2}\omega ^{2}+r_{0}=0
\end{equation}%
where

$r_{6}=n_{43}^{2}-2n_{42},$ $%
r_{4}=n_{42}^{2}+2n_{40}-2n_{43}n_{41}-n_{22}^{2},$

$r_{2}=n_{41}^{2}-2n_{42}n_{40}+2n_{22}n_{20}-n_{21}^{2},$ $%
r_{0}=n_{40}^{2}-n_{20}^{2}.$

If $\omega _{0}$ is a positive root of $\left( 15\right) $ then there is a
Hopf bifurcation and the value of $\tau _{0}$ is given by:
\begin{equation}
\tau _{0}=\dfrac{1}{\omega _{0}}arctg\dfrac{a_{1}a_{4}\omega _{0}+a_{2}a_{3}%
}{-a_{1}a_{3}+a_{2}a_{4}\omega _{0}},
\end{equation}%
where $a_{1}=\omega _{0}^{4}-n_{42}\omega _{0}^{2}+n_{40},$ $%
a_{2}=-n_{43}\omega _{0}^{3}+n_{41}\omega _{0},$
$a_{3}=n_{22}\omega _{0}^{2}-n_{20},$ $a_{4}=n_{21}\omega _{0}.$

We can conclude with the following theorem:

\begin{theorem}(i) If $\omega _{0}$ is a positive root of $\left( 15\right) $ and $Re
(\frac{d\lambda }{d\tau })_{\lambda =i\omega _{0},\tau =\tau _{0}}\neq 0,$ where $%
\tau _{0}$ is given by $\left( 16\right) $, then a Hopf
bifurcation occurs at the stationary state $\left( x_{1}^{\ast
},x_{2}^{\ast },z_{1}^{\ast },z_{2}^{\ast }\right) $ as $\tau $
passes through $\tau _{0}.$

(ii) If conditions (11) hold and $n_0>0$, then the stationary
state is asymptotically stable for any $\tau>0$.
\end{theorem}

\section{Numerical simulation}

\hspace{0.4cm} For the numerical simulation we use Maple 11 and
the following data: $q=0.3$, $s=40$, $t_1=0.16$, $c_1=0.2$,
$c_2=2$, $k_1=0.05$, $k_2=0.01$, $h_1=0.05$, $h_2=0.01$. The
stationary state is: $x_{1}^{\ast }=0.34710$, $x_{2}^{\ast
}=0.0347$, $z_{1}^{\ast }=0.85075$, $z_{2}^{\ast }=0.03257$.

For $\tau=0$ the Routh-Hurwicz conditions are satisfied. Then, the
stationary state is stable.

The positive solution of (15) is $\omega _{0}=0.010083$ and
$\tau_0=164.5979$. For $\tau\in(0,\tau_0)$ the stationary state is
asymptotically stable and for $\tau\in[\tau_0,\infty)$ the
stationary state is unstable. For $\tau=\tau_0$ there is a Hopf
bifurcation.

For $q=0.3$, $s=40$, $t_1=0.16$, $c_1=0.2$, $c_2=1.5$, $k_1=0.05$,
$k_2=0.01$, $h_1=0.05$, $h_2=0.01$. The stationary state is:
$x_{1}^{\ast }=0.4359$, $x_{2}^{\ast }=0.05813$, $z_{1}^{\ast
}=0.824019$, $z_{2}^{\ast }=0.059313$.

For $\tau=0$ the Routh-Hurwicz conditions are satisfied. Then, the
stationary state is stable.

The equation (15) has no positive root. Then, the stationary state
is asymptotically stable for any $\tau>0$.

\section{Conclusions}

\hspace{0.4cm} In the static model with tax evasion, the
parameters $q$ and $s$ characterize the behavior of the firms with
respect to evasion. The presented figures allow the analysis of
the declared revenues and the profits with respect to $s$.

For the dynamic model with tax evasion, using Routh-Hurwitz
criterion we have determined the conditions for which the
stationary state is asymptotically stable$.$

For the dynamic model with tax evasion and time delay, using the delay $\tau
$ as a bifurcation parameter we have shown that a Hopf bifurcation occurs
when $\tau $ passes through a critical value $\tau _{0}.$

The direction of the Hopf bifurcation, the stability and the period of the
bifurcating periodic solutions will be analyzed in a future paper.

The findings of the present paper can be extended in an oligopoly case.

\section{References}

\hspace{0.4cm} [1]. C. Chiarella, F. Szidarovszky, \textit{\ On
the asymptotic behavior of dynamic rent-seeking games}, electronic
Journal: Southwest Journal of Pure and Applied Mathematics, Issue
2:17-27, 2000.

[2]. A. Gebauer, C.W. Nam, R. Parsche,\textit{\ VAT evasion and
its consequence for macroeconomic clearing in the EU},
FinanzArchiv, 2005.

[3]. L. Goerke, M. Runkel, \textit{Indirect Tax Evasion and
Globalization,} FinanzArchiv, 2006.

[4]. R.H. Gordon, S.B. Nielsen, \textit{Tax evasion in an open
economy: value-added vs. income taxation,} Journal of Public
Economics 66, 173-197, 1997.

[5]. K. Matthews, J. Lloyd-Williams, \textit{The VAT-evading firm
and VATevasion: an empirical analysis,} International Journal of
the Economics of Business, 8, 39-49, 2001.

[6]. G. Mircea, M. Neam\c{t}u, D. Opri\c{s}, \textit{Hopf
bifuraction for dynamical system with time delay and application,} Ed. Mirton, Timi\c{s}%
oara, 2004.

[7]. M. Neam\c{t}u, C. Chil\u{a}rescu, \textit{On the Asymtotic
Behavior of Dynamic Rent-Seeking Games with Distributed Time},
Proceedings of the Sixth International Conference on Economic
Informatics, Bucharest, 8-11 may, Ed. Economica, 216-223, 2003.

[8]. M. Neam\c{t}u, D. Opri\c{s}, \textit{Rent Seeking Dynamic
Games with Delay}, "Globalizarea \c{s}i educa\c{t}ia economic\u{a} universitar\u{a}" Ia%
\c{s}i, 207-212, Editura Sedcom Libris, 2002.

[9]. C. Roc\c{s}oreanu, \textit{The bifurcations of continues
dynamical systems. Applications in economy and biology,} Ed.
Universitaria Craiova, 2006.

[10]. L. Xu, F. Szidarovszky, \textit{The Stability of Dynamic
Rent-Seeking Games}, International Game Theory Reviews, No 1:
87-102, 1999.

\end{document}